%
%
%

\def\LaTeX{%
  \let\Begin\begin
  \let\End\end
  \let\salta\relax
  \let\finqui\relax
  \let\futuro\relax}

\def\UK{\def\our{our}\let\sz s}
\def\USA{\def\our{or}\let\sz z}



\LaTeX

\USA


\salta

\documentclass[twoside,12pt]{article}
\setlength{\textheight}{24cm}
\setlength{\textwidth}{16cm}
\setlength{\oddsidemargin}{2mm}
\setlength{\evensidemargin}{2mm}
\setlength{\topmargin}{-15mm}
\parskip2mm


\usepackage{amsmath}
\usepackage{amsthm}
\usepackage{amssymb}
\usepackage[mathcal]{euscript}

\usepackage[usenames,dvipsnames]{color}
%
%
%
\def\gianni{\color{green}}
\def\pier{\color{red}}
\def\elena{\color{blue}}
\def\elenabis{\color{cyan}}
%
%
\let\gianni\relax
\let\pier\relax
\let\elena\relax
\let\elenabis\relax

\bibliographystyle{plain}


%

\finqui

\def\Beq{\Begin{equation}}
\def\Eeq{\End{equation}}
\def\Bsist{\Begin{eqnarray}}
\def\Esist{\End{eqnarray}}

\def\Bthm{\Begin{theorem}}
\def\Ethm{\End{theorem}}

\def\Bprop{\Begin{proposition}}
\def\Eprop{\End{proposition}}
\def\Bcor{\Begin{corollary}}
\def\Ecor{\End{corollary}}
\def\Brem{\Begin{remark}\rm}
\def\Erem{\End{remark}}

\def\Bdim{\Begin{proof}}
\def\Edim{\End{proof}}
\let\non\nonumber




\def\step #1 \par{\medskip\noindent{\bf #1.}\quad}

\def\Step #1 \par{\bigskip\leftline{\bf #1}\nobreak\medskip}


\def\Lip{Lip\-schitz}
\def\holder{H\"older}
\def\aand{\quad\hbox{and}\quad}

\def\lsc{lower semicontinuous}

\def\lhs{left-hand side}
\def\rhs{right-hand side}

\def\wk{well-known}


\def\generaliz{generali\sz}

\def\organiz{organi\sz}

\def\bhv{behavi\our}


\def\multibold #1{\def\arg{#1}%
  \ifx\arg\pto \let\next\relax
  \else
  \def\next{\expandafter
    \def\csname #1#1#1\endcsname{{\bf #1}}%
    \multibold}%
  \fi \next}

\def\pto{.}

\def\multical #1{\def\arg{#1}%
  \ifx\arg\pto \let\next\relax
  \else
  \def\next{\expandafter
    \def\csname cal#1\endcsname{{\cal #1}}%
    \multical}%
  \fi \next}


\def\multimathop #1 {\def\arg{#1}%
  \ifx\arg\pto \let\next\relax
  \else
  \def\next{\expandafter
    \def\csname #1\endcsname{\mathop{\rm #1}\nolimits}%
    \multimathop}%
  \fi \next}

\multibold
qwertyuiopasdfghjklzxcvbnmQWERTYUIOPASDFGHJKLZXCVBNM.

\multical
QWERTYUIOPASDFGHJKLZXCVBNM.

\multimathop
dist div dom meas sign supp .


\def\accorpa #1#2{\eqref{#1}--\eqref{#2}}
\def\Accorpa #1#2 #3 {\gdef #1{\eqref{#2}--\eqref{#3}}%
  \wlog{}\wlog{\string #1 -> #2 - #3}\wlog{}}


\def\tonde #1{\left(#1\right)}

\def\graffe #1{\mathopen\{#1\mathclose\}}

\def\<#1>{\mathopen\langle #1\mathclose\rangle}
\def\norma #1{\mathopen \| #1\mathclose \|}

\def\normaV #1{\norma{#1}_V}
\def\normaH #1{\norma{#1}_H}

\let\meno\setminus

\def\iot {\int_0^t}
\def\ioT {\int_0^T}
\def\iO{\int_\Omega}
\def\intQt{\int_{Q_t}}
\def\intQ{\int_Q}

\def\dt{\partial_t}
\def\dn{\partial_\nu}
\def\ds{\,ds}

\def\cpto{\,\cdot\,}

\def\checkmmode #1{\relax\ifmmode\hbox{#1}\else{#1}\fi}
\def\aeO{\checkmmode{a.e.\ in~$\Omega$}}
\def\aeQ{\checkmmode{a.e.\ in~$Q$}}

\def\aat{\checkmmode{for a.a.~$t\in(0,T)$}}


\def\erre{{\mathbb{R}}}




\def\genspazio #1#2#3#4#5{#1^{#2}(#5,#4;#3)}
\def\spazio #1#2#3{\genspazio {#1}{#2}{#3}T0}

\def\L {\spazio L}
\def\H {\spazio H}
\def\W {\spazio W}

\def\C #1#2{C^{#1}([0,T];#2)}

\def\Vp{V_0^*}
\def\Vz{V_0}


\def\Lx #1{L^{#1}(\Omega)}
\def\Hx #1{H^{#1}(\Omega)}
\def\Wx #1{W^{#1}(\Omega)}
\def\Wxz #1{W^{#1}_0(\Omega)}
\def\Cx #1{C^{#1}(\overline\Omega)}

\def\HxG #1{H^{#1}(\Gamma)}

\def\Luno{\Lx 1}
\def\Ldue{\Lx 2}
\def\Linfty{\Lx\infty}
\def\Lq{\Lx4}
\def\Huno{\Hx 1}
\def\Hdue{\Hx 2}
\def\Hunoz{{H^1_0(\Omega)}}


\def\LQ #1{L^{#1}(Q)}


\let\theta\vartheta
\let\epsilon\varepsilon

\let\TeXchi\chi                         
\newbox\chibox
\setbox0 \hbox{\mathsurround0pt $\TeXchi$}
\setbox\chibox \hbox{\raise\dp0 \box 0 }
\def\chi{\copy\chibox}


\let\eps\tau
\let\ka\kappa


\def\thetae{\theta_\eps}
\def\chie{\chi_\eps}
\def\xie{\xi_\eps}
\def\thetaz{\theta_0}
\def\chiz{\chi_0}
\def\thetaG{\theta_\Gamma}
\def\thetaH{\theta_\calH}
\def\ue{u_\eps}
\def\thetamin{\theta_*}
\def\thetamax{\theta^*}
\def\phie{\phi_\eps}

\def\Beta{\widehat{\vphantom t\smash\beta\mskip2mu}\mskip-1mu}

\def\kaz{\ka_0}
\def\kalim{\ka_0'}
\def\ke{k_\eps}
\def\ketriv{\smash{\tilde k_\eps}}
\def\kamin{\ka_*}
\def\kamax{\ka^*}
\def\Ln{\mathop {\rm Ln}}
\def\chart{{\cal C}_t}

\def\todx{\mathrel{\scriptstyle\searrow}}

\Begin{document}


\title{\bf Singular limit of an integrodifferential system 
related to the entropy balance\footnote{{\pier {\bf 
Acknowledgments.}\quad\rm The authors gratefully 
acknowledge the financial support of the MIUR-PRIN Grant 
2010A2TFX2 ``Calculus of variations'' and of the IMATI of CNR in Pavia.}}}
\author{}
\date{}
\maketitle
\begin{center}
\vskip-2cm
{\large\bf Elena Bonetti$^{(1)}$}\\
{\normalsize e-mail: {\tt elena.bonetti@unipv.it}}\\[.25cm]
{\large\bf Pierluigi Colli$^{(1)}$}\\
{\normalsize e-mail: {\tt pierluigi.colli@unipv.it}}\\[.25cm]
{\large\bf Gianni Gilardi$^{(1)}$}\\
{\normalsize e-mail: {\tt gianni.gilardi@unipv.it}}\\[.45cm]
$^{(1)}$
{\small Dipartimento di Matematica ``F. Casorati'', Universit\`a di Pavia}\\
{\small via Ferrata 1, 27100 Pavia, Italy}\\[.8cm]
{\large\it Dedicated to Professor Mauro Fabrizio,\\[.5mm]
nice friend and valuable collaborator,\\[1.5mm]
wishing him the very best ``ad multos annos''}\\[.8cm]
\end{center}


\Begin{abstract}{\pier A thermodynamic model describing phase transitions with thermal memory, {\elenabis  in terms} of an entropy equation and a momentum balance for the  
microforces, is adressed. Convergence results and error estimates are proved  
for the related integrodifferential system of {\elenabis PDE} as the sequence of memory kernels converges  to a multiple of a Dirac delta, {\elenabis in a suitable sense}.} 
\\[2mm]
\End{abstract}


\salta

\pagestyle{myheadings}
\newcommand\testopari{\sc Bonetti \ --- \ Colli \ --- \ Gilardi}
\newcommand\testodispari{\sc  Singular limit of an integrodifferential system}
\markboth{\testodispari}{\testopari}

\finqui

%
%
\section{Introduction}
\label{Intro}
\setcounter{equation}{0}
In this paper, {\elena we deal with the following integrodifferential  PDE system, describing a  
phase transition process in the case when thermal
memory effects are included. Indeed, here 
$\thetae$ stands for the absolute temperature, 
$\chie$ for the phase parameter, and $\ke$ for a memory kernel}
\Bsist
  && \dt \bigl( \ln\thetae + \lambda(\chie) \bigr)
  - \kaz \Delta\thetae - \Delta (\ke*\thetae) = f
  \label{Iprimaeps}
  \\
  && {\elena \dt\chie - \Delta\chie + \xi_\eps + \sigma'(\chie) = \lambda'(\chie) \, \thetae,\quad\xi_\eps\in\beta(\chie)}
  \label{Isecondaeps}
  \\
  && \thetae|_\Gamma = \thetaG
  \aand
  \dn\chie|_\Gamma = 0
  \label{Ibceps}
  \\
  && \ln\thetae(0) = \ln\thetaz
  \aand
  \chie(0) = \chiz .
  \label{Icauchyeps}
\Esist
\Accorpa\Ipbleps Iprimaeps Icauchyeps
Each of the partial differential equations \accorpa{Iprimaeps}{Isecondaeps} 
is meant to hold in a three-dimensional bounded domain $\Omega$, 
endowed with a smooth boundary~$\Gamma$, and in some time interval~$(0,T)$.
{\pier In~\eqref{Iprimaeps}, the memory kernel~$\ke$ may depend 
on a positive parameter $\eps$.}
Moreover, the symbol~$*$ denotes the usual time convolution formally defined by
$(a*b)(t)=\iot a(t-s)\,b(s)\ds$ for functions that depend just on time, 
and then extended to functions that also depend on space.
Furthermore, $f$~is some given source term.
In~\eqref{Isecondaeps}, $\beta$~is a maximal monotone graph in~$\erre^2$,
while $\lambda$ and $\sigma$ are real functions defined on the whole of~$\erre$.
{\pier The boundary conditions \eqref{Ibceps} must be satisfied in $\Gamma 
\times (0,T)$, while the initial conditions \eqref{Icauchyeps} are written for the functions $\ln \thetae$ and $\chie$: of course,} $\thetaG$, $\thetaz$, and $\chiz$ are given boundary and initial data.
 
{\elena Equation \eqref{Iprimaeps} may be interpreted as an entropy balance equation. Note in particular
that the equation is singular with respect to the temperature, mainly for the presence of the logarithm, forcing the temperature
to assume only positive values (which is in accordance with physical consistency). Similar systems have been studied in the literature
from the point of view of the existence and regularity of solutions (see, among the others, \cite{BCF,BCFG2,BCFG1,BCFG3,BFR,GR}).}

{\pier The well-posedness of a proper variational formulation of \Ipbleps\ has been
proved in  \cite{BCFG1}. Here,} our main goal is the following.
By assuming that $\ke$ 
converges to $\kalim\delta$ at $\eps\todx0$ in a suitable sense,
where $\delta$ is the Dirac mass at the origine of the real line and $\kalim$ is a real constant
satisfying $\ka:=\kaz+\kalim>0$,
we prove that the solution $(\thetae,\chie)$ to \Ipbleps\ converges in a proper topology
to the solution $(\theta,\chi)$ of the problem stated below
\Bsist
  && \dt \bigl( \ln\theta + \lambda(\chi) \bigr)
  - \ka \Delta\theta = f
  \label{Iprima}
  \\
  && {\elena \dt\chi - \Delta\chi + \xi + \sigma'(\chi) = \lambda'(\chi) \, \theta,\quad \xi\in\beta(\chi)}
  \label{Iseconda}
  \\
  && \theta|_\Gamma = \thetaG
  \aand
  \dn\chi|_\Gamma = 0
  \label{Ibc}
  \\
  && \ln\theta(0) = \ln\thetaz
  \aand
  \chi(0) = \chiz \,.
  \label{Icauchy}
\Esist
\Accorpa\Ipbl Iprima Icauchy
{\elena {\pier This convergence result is obtained by the} use of an a priori estimates technique and passage to the limit arguments, based on monotonicity and compactness}. Moreover, {\pier an {\it error estimate}, 
i.e. an estimate of suitable norms or quantities involving 
the difference of solutions,} is shown.

Our paper is \organiz ed as follows.
In the next section, we {\pier discuss a derivation of the system \accorpa{Iprima}{Iseconda} from the basic laws of thermomechanics.}
Section~\futuro\ref{Statement} is devoted to
the statement of our assumptions and of our results
on the mathematical problem.
In Section~\futuro\ref{Auxiliary}, we present some auxiliary material
that is needed for the proof of our convergence Theorem~\futuro\ref{Convergenza}, mainly.
The last section is devoted to the proofs of
the above theorem and of the error estimate {\elenabis stated in Theorem \ref{Errore}}.


\section{The model}
\label{Model}
\setcounter{equation}{0}

\newcommand\flussond{-{\bf Q}^{nd}}
\newcommand\flussod{-{\bf Q}^d}

\newcommand\storiee{{\cal S}_\tau}
\newcommand\he{h_\tau} 
\newcommand\histem{\widetilde{\nabla\theta}^t}

In this section, we  briefly {\elena introduce} the modeling derivation of the equations  \accorpa{Iprimaeps}{Isecondaeps} 
and discuss {\pier the} convergence to \accorpa{Iprima}{Iseconda}, {\elena as the parameter $\tau$ (in the memory kernel) {\pier tends to $0$}. Here,
the argument is mainly developed from a physical point of view, while we refer to subsequent sections for a more {\pier precise setting}
of analytical assumptions and comments}. In particular, 
we aim to focus on the fact that \eqref{Iprimaeps} accounts for thermal evolution involving memory effects, 
on the basis of the memory kernel $\ke$. 

Materials with thermal memory have been deeply studied in the literature, both from a modeling 
and analytical point of view. We refer, in particular, to the 
approach by Gurtin and 
Pipkin (see \cite{GPmodel}) for thermal memory materials. 
{\elena Several authors have investigated phase transitions in special materials with thermal memory, both concerning modeling and analysis.
{\elenabis For a fairly complete and detailed presentation of this kind of problems, let us mention the very recent monograph \cite{librofabrizio}.} Now, we  combine} {\pier thermal memory} 
with a new theory for phase transitions models, based on a generalization of the principle of virtual powers (see \cite{Fremond}).
The idea is that  micro-forces, which are responsible for the phase transition, have to be included in the whole energy
balance of the system. Consequently, the phase (evolution) equation is derived as a  micro-forces balance equation and
it is {\pier coupled} with an entropy evolution equation.
This approach has been recently   {\elena investigated in the literature by several authors} 
(among the others, we mainly refer to the papers~\cite{BCF} and~\cite{BFR}, in which the derivation of the model is detailed in the case 
when possible thermal memory effects are included, as in equations~\accorpa{Iprimaeps}{Isecondaeps}).   

{\elena Indeed,} let us recall that in \cite{BCF} the theory by Gurtin-Pipkin 
is considered, allowing the free energy functional 
to depend on the past history of the temperature gradient. The resulting  functional accounts for   
non-dissipative contributions in the heat flux, which may be combined with additional dissipative 
instantaneous contributions coming from a pseudo-potential of dissipation. The use of an entropy
balance {\elena has been recovered, in this approach,} from a rescaling (with respect to the absolute temperature) of the energy balance, 
under the small perturbations assumption (see also \cite{BCFG2, BCFG1}).
In \cite{BFR} a fairly general theory is introduced. The model is derived by  a dual approach (mainly in the sense of convex analysis) 
in which the entropy 
and the history of the entropy flux are chosen as state variables 
(together with the phase parameter and possibly its gradient). 
Then, the dissipative functional  is written in terms of a dissipative contribution in the entropy flux and for the 
time derivative of the phase parameter. 

Let us point out that the above mentioned  approach is not far from the theory proposed by Green-Naghdi \cite{GreenNaghdi}  
and Podio-Guidugli~\cite{Podio1},  in which some 
{\it thermal displacement} is introduced as state variable  (it is a primitive of the temperature) and 
the equations come from a  generalization of the principle of virtual powers, in which thermal 
forces are included. 
As a consequence, in this framework, the entropy equation is  formally obtained as a momentum balance (i.e., a balance of thermal forces acting in the system). {\pier The reader may also examine~\cite{CC1, CC2}, where some asymptotic analyses are carried out to find the interconnections among peculiar Green and Naghdi types.} 

We aim {\pier to} observe that the model we are investigating actually   may {\elena be obtained by combining} the above {\elena two} 
theories, {\elena i.e.} generalizing  
the principle of virtual powers accounting for  microforces 
as well as thermal stresses. {\elena Let us present our position.} First, 
{\pier we} specify the expression of the power of internal forces.
The power of interior forces is written for any virtual 
micro-velocity $\gamma$ and thermal velocity $v$,
as follows
\Beq\label{Pi}
{\cal P}_i=\int_\Omega B\gamma+{\bf H}\cdot\nabla\gamma+{\bf Q}\cdot\nabla v,
\Eeq
where $B$ and ${\bf H}$ are interior forces responsible for the phase transition {\elenabis (as introduced in \cite{Fremond})}, and ${\bf Q}$ stands for a thermal stress 
(corresponding to the entropy flux {\elenabis by} \cite{Podio1}).  
Hence, the resulting balance equations are  written as momentum balance equations. 
It is assumed that an external {\elena (density of)} entropy  source $f$ is applied. A {\it thermal momentum} is introduced to measure  reluctance to 
{\pier the} order of the system (in analogy with the mechanical momentum measuring reluctance to quiet).
We prescribe that it is given by the entropy $s$.
It results that (see~\eqref{Pi})
\Beq\label{entropy}
s_t+\hbox{div }{\bf Q}=f.
\Eeq
As far as the microscopic momentum balance is concerned, we assume that no acceleration and no external force are contributing, so that we have
\Beq\label{momentum}
B-\hbox{div }{\bf H}=0.
\Eeq
Henceforth, \accorpa{entropy}{momentum}  are
combined with suitable boundary conditions. {\elenabis As usual, we assume that the flux through the boundary  ${\bf H}\cdot{\bf n}$  is null, while (mainly for analytical reasons) we prescribe a known temperature 
on the boundary.} 

The entropy $s$, 
the entropy flux ${\bf Q}$, and the new interior forces $B$ and ${\bf H}$ are recovered by suitable energy and dissipation 
functionals,
that we are going to make precise, in terms of state variables. The state variables are related to 
the equilibrium of the thermodynamical system: they are
the absolute temperature $\theta$, the phase parameter $\chi$, 
the gradient $\nabla\chi$ {\elena (actually accounting for local interactions)}, and the history variable $\histem$, which is defined as
\Beq
\histem(s)=\int_{t-s}^t\nabla\theta(r)dr{\pier , \quad s>0 .} 
\Eeq
As in \cite{BCF}, we assume that the free energy of the system 
(depending on $(\theta,\histem,\chi,\nabla\chi)$) is {\pier split} into two contributions: 
the first is related to present variables at time $t$ ($\Psi_P$), 
the second accounts for some history in the system ($\Psi_H$), measured through a memory kernel (related to $\ke$ in {\elena the} equations). 
In particular, the history contribution of the free energy  is given by
\Beq\label{psiH}
\Psi_H(\histem)=\frac 1 2|\histem|^2_{\storiee}
\Eeq
where $\storiee$ is the space of the past histories  
(as it is introduced in the theory of thermal memory materials by Gurtin and Pipkin), defined by
\Beq
\storiee:=\{{\bf f}:(0,+\infty)\rightarrow\erre^3\hbox{ measurable s.t. }\int_0^{+\infty}\he(s)|{\bf f}(s|^2 ds<+\infty\} .
\Eeq
Here, $\he:(0,+\infty)\rightarrow(0,+\infty)$ (possibly depending on a parameter $\tau$) is a continuous, decreasing function 
such that 
\Beq
\int_0^{+\infty}s^2 \he(s)ds<+\infty.
\Eeq
The space $\storiee$ is endowed with the natural norm
\Beq\label{normaS}
|{\bf f}|^2_{\storiee}=\int_0^{+\infty} \he 
(s)|{\bf f}(s)|^2 ds
\Eeq
and the related scalar product is $({\bf v},{\bf u})_{\storiee}=\int_0^{+\infty} \he (s){\bf v}(s)\cdot{\bf u}(s)ds$.
Let us comment that in our system, to derive \eqref{Iprimaeps},  we have introduced  a kernel  $\ke$ such that $-\ke'=\he$.
More precisely, let $\ke : (0,+\infty)\rightarrow \erre$
and require that
\Beq
\label{ipoke}
\ke\in W^{2,1} (0,+\infty) , \quad \lim_{s\rightarrow + \infty}
\ke(s)=0.
\Eeq
Hence, by virtue of the assumptions on $\he$ we also have
\Beq
\ke'\leq0\quad\hbox{and}\quad \ke''\geq0\quad\hbox{a.e. in }\, (0,+\infty).
\Eeq
Note that $\ke'(t)$ vanishes for $t$ going to $+\infty$  
and that
$\ke$ is a non-increasing function with $\ke(0)\geq 0,$ and
in the case $\ke(0)=0$ one has $\ke \equiv 0$. These assumptions on $\ke$ actually ensure that the model is thermodynamically consistent,
as it is detailed in \cite{BCF}.


Then, the free energy functional $\Psi_P$ (written at the present time $t$) is addressed
\Beq
\Psi_P(\theta,\chi,\nabla\chi)=c_V\theta(1-\ln\theta)-\lambda(\chi)\theta+\sigma(\chi)+\widehat\beta(\chi)+\frac 1 2|\nabla\chi|^2 
\Eeq
where $c_V>0$ (in the sequel let us take $c_V=1$) is the specific heat, $\sigma$ and $\lambda$ are sufficiently smooth 
functions (with $\lambda'(\chi)$ denoting the latent heat), 
$\widehat\beta$ is a proper convex and lower semicontinuous function, possibly accounting for internal 
constraints on the phase variable~$\chi$. {\elena For instance, a fairly classical choice is $\widehat\beta(\chi)=I_{[0,1]}(\chi)${\pier , which is equal to $0$ if $\chi\in[0,1]$
and takes value $+\infty$ elsewhere} ({\pier thus} forcing $\chi\in[0,1]$).}

Dissipation is rendered in terms of the time derivative $\chi_t$ and of the dissipative variable~$\nabla\theta$. 
It is derived by a pseudo-potential of dissipation (in the sense of Moreau, i.e. a convex, non-negative function assuming its minimum  0
for null dissipation):
\Beq\label{pseudovec}
\Phi(\chi_t,\nabla\theta)=\frac 1 2|\chi_t|^2+\frac \kaz 2 |\nabla\theta|^2. 
\Eeq
{\pier Note that, in order to ensure the validity of the second principle of thermodynamics,} it is required that $\kaz\geq0$.

Now, we are in a position to recover our system, after specifying constitutive relations for the involved physical quantities.
We have that 
\Beq\label{defs}
s=-\frac{\partial\Psi}{\partial\theta}=\ln\theta+\lambda(\chi)
\Eeq
and 
\Beq
B{\elena{}\in{}}\frac{\partial\Psi}{\partial\chi}+\frac{\partial\Phi}{\partial\chi_t}=\beta(\chi)+\sigma'(\chi)-\lambda'(\chi)\theta {\pier {}+ \chi_t}
\Eeq
$\beta$ being the subdifferential (in the sense of convex analysis) 
{\pier of $\widehat\beta$}, and 
\Beq
{\bf H }=\frac{\partial\Psi}{\partial(\nabla\chi)}=\nabla\chi.
\Eeq
Hence, the entropy flux vector ${\bf Q}$ is specified by 
\Beq\label{defQ}
-{\bf Q}=(\flussond)+(\flussod)
\Eeq
where  $\flussond$ results to be defined in $\storiee$. It is obtained {\pier
taking the derivative in $\storiee$ of} the history functional 
with respect to the history variable{\pier . Integrating by parts in time, 
using the Fr\'echet derivative, and  exploiting 
 the hypotheses on $\ke$ (see \cite{BCF} for any further detail) lead to}
\Beq\label{flussnd}
\flussond {\pier (t)}  =-\int_{-\infty}^t\ke(t-s)\nabla\theta(s)ds.
\Eeq
{\pier We have now} to make precise the dissipative part of the entropy flux 
\Beq\label{defqduno}
\flussod= \frac{\partial\Phi}{\partial\nabla\theta}=\kaz\nabla\theta.
\Eeq

Hence, equations \accorpa{Iprimaeps}{Isecondaeps} 
are obtained by \eqref{entropy} and \eqref{momentum} {\elena exploiting} the above introduced 
constitutive relations. We point out that in \eqref{Iprimaeps}, the past history contribution of \eqref{flussnd} 
(actually its divergence), 
i.e. $\int_{-\infty}^0\ke(t-s)\nabla\theta(s)ds$, is assumed to be known and included in the external entropy source $f$ 
(we have used the same notation as in 
\eqref{entropy} for the sake of simplicity).

As we have already pointed out in the Introduction, the main aim of this paper is  to investigate the {\elena asymptotic} {\gianni\bhv} of 
system  \accorpa{Iprimaeps}{Isecondaeps} as
the thermal memory kernel converges to $\kalim\delta$, 
$\delta$ being the Dirac mass at the origin of the real line  and $\kalim>0$, in a suitable sense
\Beq
\ke\rightarrow\kalim\delta.
\Eeq

{\elenabis More precisely, we are interested in proving that solutions to the system \accorpa{Iprimaeps}{Isecondaeps} converge to solutions to \accorpa{Iprima}{Iseconda} (at least in some weak topology). Let us briefly comment that the system \accorpa{Iprima}{Iseconda}, obtained in our proof as a suitable limit of \accorpa{Iprimaeps}{Isecondaeps}, can be actually derived by an analogous procedure as {\pier the one} we have performed 
to formally derive \accorpa{Iprimaeps}{Isecondaeps}. Indeed, \accorpa{Iprima}{Iseconda} follow from \eqref{entropy} and \eqref{momentum} {\pier when} exploiting \eqref{defs}-\eqref{defQ}{\pier . Here,} the new energy and dissipative functionals are $\Psi=\Psi_P$ (i.e., no history contribution 
of type \eqref{psiH} in the free energy is given) and (cf. \eqref{pseudovec})
\Beq\label{pseudonew}
\Phi(\chi_t,\nabla\theta)=\frac 1 2|\chi_t|^2+\frac {(\kaz+\kalim)} 2 |\nabla\theta|^2.
\Eeq
 In particular, it results that in \eqref{defQ} ${\bf Q}^{nd}={\bf 0}$, while (due to \eqref{defqduno}) $-{\bf Q}^d=(\kaz+\kalim)\nabla\theta$.}


\section{Statement of the mathematical problem}
\label{Statement}
\setcounter{equation}{0}

In this section, we make our assumptions precise and state our results.
First of all,
we assume $\Omega$ to be a bounded connected open set in~$\erre^3$
(lower-dimensional cases could be considered with minor changes)
whose boundary~$\Gamma$ is supposed to be smooth.
Next, we fix a final time $T\in(0,+\infty)$
and set:
\Bsist
  && Q := \Omega\times(0,T) , \quad
  \Sigma := \Gamma\times(0,T)
  \label{defQS}
  \\
  && V := \Huno, 
  \quad \Vz :=\Hunoz ,
  \quad H := \Ldue 
  \label{defVH}
  \\
  && W := \graffe{v\in\Hdue:\ \dn v = 0 \ \hbox{on $\Gamma$}},
  \label{defW}
\Esist
{\elena $\dn$ denoting the normal derivative operator on the boundary}.
\Accorpa\Defspazi defVH defW
We endow the spaces \Defspazi\ with their standard norms,
for which we use a self-explained notation like $\normaV\cpto$.
Moreover, for $p\in[1,+\infty]$, we write $\norma\cpto_p$ for the usual norm
in~$L^p(\Omega)$; as no confusion can arise, the symbol $\norma\cpto_p$
is used for the norm in $L^p(Q)$ as well.
In the sequel, the same symbols are used for powers of the above spaces {\elena and the corresponding natural induced norms}.
It is understood that $H\subset\Vp$ as usual, i.e.,
any element $u\in H$ is identified
with the functional $\Vz\ni v\mapsto\iO uv$ which actually belongs 
to the dual space~$\Vp=\Hx{-1}$ of~$\Vz$.
We observe that $\L2\Vp$ coincides with the dual space of $\L2\Vz$
and use the symbol $\<\cpto,\cpto>$ for the corresponding duality pairing.

\bigskip

As far as the structure of the system is concerned {\elena (see \eqref{Iprimaeps},~\eqref{Iprima} and \eqref{Isecondaeps},~\eqref{Iseconda})}, 
we are given {\elena the} three
functions $\Beta$, $\lambda$, $\sigma$, 
{\elena the}~constant $\kaz$ and {\elena the} memory kernel $\ke$ depending on the parameter $\eps>0$
and we assume that the conditions listed below are satisfied.
\Bsist
  && \hbox{$\Beta:\erre\to[0,+\infty]$ is convex, proper, lower
      semicontinuous, and $\Beta(0)=0$}
    \qquad
  \label{hpBeta}
  \\
  && \hbox{$\lambda,\sigma\in C^1(\erre)$ and
      $\lambda',\sigma'$ are \Lip\ continuous}
  \label{hpls}
  \\
  && \kaz > 0 \aand 
  \ke \in W^{1,1}(0,T).
  \label{hpregk}
\Esist
We define the graph $\beta$ in $\erre\times\erre$ by
\Beq
  \beta := \partial\Beta
  \label{defbeta}
\Eeq
\Accorpa\Hpstruttura hpBeta defbeta
and note that $\beta$ is maximal monotone and that $\beta(0)\ni0$. 
In the sequel, we write $D(\Beta)$ and $D(\beta)$ for the effective domains
of~$\Beta$ and~$\beta$, respectively, and we use the same {\pier symbol 
$\beta$ for 
the maximal monotone operators induced on $L^2$ spaces.} 

As far as the data of our problem are concerned,
we assume that the functions $f$, $\thetaG$, $\thetaz$, $\chiz$ 
and the constants $\thetamin$ and $\thetamax$ are given such~that
\Bsist
  && f \in \L2H 
  \label{hpf}
  \\
  && \thetaG \in  {\pier \H1{\HxG{-1/2}} \cap \W{1,1}{L^\infty(\Gamma)}
    \cap \C0{\HxG{1/2}}} \qquad
  \label{regthetaG}
  \\
  && 0 < \thetamin \leq \thetamax < +\infty
  \label{hpthetaminmax}
  \\
  && \thetamin \leq \thetaG \leq \thetamax
    \quad \hbox{a.e.\ on $\Gamma\times(0,T)$}
  \label{hpthetaG}
  \\
  && \thetaz \in \Linfty, \quad
    \thetamin \leq \thetaz \leq \thetamax \quad \aeO
  \label{hpthetaz}
  \\
  && \chiz \in V \aand \Beta(\chiz) \in \Luno.
  \label{hpchiz}
\Esist
\Accorpa\Hpdati hpf hpchiz
The function $\thetaG$ is the boundary datum for the temperature and we
would like to consider a function $u:=\theta-\thetaH$ vanishing on the boundary as associated unknown function.
Hence, a~natural choice of $\thetaH$ is the harmonic extension of~$\thetaG$, so that $\Delta u=\Delta\theta$.
Therefore, we define $\thetaH:Q\to\erre$ as~follows
\Beq
  \thetaH(t) \in V, \quad \Delta \thetaH(t) = 0 ,
  \aand \thetaH(t)|_\Gamma = \thetaG(t) \quad \aat.
  \label{defthetaH}
\Eeq
{\elena For the regularity of $\thetaH$ (induced by \eqref{regthetaG}) see {\pier the} subsequent Proposition \ref{StimethetaH}}.

Next, we list the a priori regularity conditions we require 
for any solution $(\theta,\chi,\xi)$ of either \Ipbleps\ or \Ipbl.
We ask~that
\Bsist
  && \theta \in \L2V \aand u := \theta - \thetaH \in \L2\Vz
  \label{regtheta}
  \\
  && \theta > 0 \quad \aeQ \aand \ln\theta \in \H1\Vp \cap \L\infty H
  \label{regln}
  \\
  && \chi \in \H1H \cap \L\infty V \cap \L2W
  \label{regchi}
  \\
  && \xi \in \L2H.
  \label{regxi}
\Esist
\Accorpa\Regsoluz regtheta regxi

At this point, we are ready to state the problems we are dealing with in a precise form.
For fixed $\eps>0$, we look for a triplet $(\thetae,\chie,\xie)$ 
satisfying \Regsoluz\ and the following system
\Bsist
  && \dt \bigl( \ln\thetae(t) + \lambda(\chie(t)) \bigr)
  - \kaz \Delta\thetae(t) - \Delta (\ke*\thetae)(t) = f(t)
  \qquad
  \non
  \\ 
  && {\pier \hskip6cm \hbox{in $\Vp$, \aat}}
  \label{primaeps}
  \\
  && \dt\chie - \Delta\chie + \xie + \sigma'(\chie) 
  = \lambda'(\chie) \, \thetae
  \aand \xie \in \beta(\chie)
  \quad \aeQ
  \qquad
  \label{secondaeps}
  \\
  && (\ln\thetae)(0) = \ln\thetaz
  \aand
  \chie(0) = \chiz .
  \label{cauchyeps}
\Esist
\Accorpa\Pbleps primaeps cauchyeps
We note that the boundary conditions \eqref{Ibceps}
are contained in \eqref{regtheta} and~\eqref{regchi}
(see the definitions \eqref{defthetaH} and \Defspazi). 
We also remark that \eqref{regln} implies that $\ln\theta$
is a continuous $\Vp$-valued function
(while {\pier no continuity of $\theta$ is known}),
so~that the Cauchy condition for $\ln\thetae$ contained in \eqref{cauchyeps} makes sense.
Similar remarks hold for the limit problem we are going to state 
(i.e., \Ipbl\ in a precise form). 

The following well-posedness result deals with a fixed $\eps>0$ and
essentially follows from~\cite{BCFG1}.
Just the notation is different, indeed.

\Bthm
\label{Buonapositura}
Assume that both \Hpstruttura\ and {\gianni\Hpdati}\ hold.
Then, there exists a unique triplet $(\thetae,\chie,\xie)$ 
satisfying \Regsoluz\ and solving problem~\Pbleps.
\Ethm

Our aim is to study the limit of the solution $(\thetae,\chie,\xie)$
as $\eps$ tends to zero, under suitable assumptions on the \bhv\ 
of the memory kernel~$\ke$.
Namely, we assume~that
\Beq
  1*\ke \to \kalim
  \quad \hbox{strongly in $L^1(0,T)$ as $\eps\todx 0$}
  \label{hpconvk}
\Eeq
for some real constant $\kalim$ 
and set
\Beq
  \ka := \kaz + \kalim \,.
  \label{defka}
\Eeq
In \eqref{hpconvk} and later on, 
we use the same symbol for any real constant
(like $1$ and~$\kalim$) 
and for the corresponding constant function.
We advice the reader that $\ka>0$ in the sequel,
either by assumption or as a consequence of some condition we require,
so that the limit problem we are going to state is parabolic with respect to~$\theta$.
Such a problem consists in looking for a triplet $(\theta,\chi,\xi)$ 
satisfying \Regsoluz\ and the following system
\Bsist
  && \dt \bigl( \ln\theta(t) + \lambda(\chi(t)) \bigr)
  - \ka \Delta\theta(t) = f(t)
  \quad \hbox{in $\Vp$, \aat}
  \qquad
  \label{prima}
  \\
  && \dt\chi - \Delta\chi + \xi + \sigma'(\chi)
  = \lambda'(\chi) \, \theta
  \aand
  \xi \in \beta(\chi) 
  \quad \aeQ
  \qquad
  \label{seconda}
  \\
  && (\ln\theta)(0) = \ln\thetaz
  \aand
  \chi(0) = \chiz .
  \label{cauchy}
\Esist
\Accorpa\Pbl prima cauchy
By just taking $\ke=0$ and replacing $\kaz$ by $\ka$ in Theorem~\ref{Buonapositura},
we obtain
\Bcor
\label{Buonaposituralimite}
Assume that both \Hpstruttura\ and {\gianni\Hpdati}\ hold
and that $\ka>0$.
Then, there exists a unique triplet $(\theta,\chi,\xi)$ 
satisfying \Regsoluz\ and solving problem~\Pbl.
\Ecor

However, we can prove a convergence result 
under further assumptions, namely
\Beq
  \ioT \bigl( \kaz v(t) + (\ke*v)(t) \bigr) v(t) \, dt
    \geq \kamin \norma v_{L^2(0,T)}^2
  \aand
  \norma\ke_{L^1(0,T)} \leq \kamax
  \qquad
  \label{hptecn}
\Eeq
for some constants $\kamin,\kamax>0$ and every $v\in L^2(0,T)$ and $\eps>0$.
By taking $v=0$ on~$(T',T)$, 
we clearly see that the time~$T$ 
can be replaced by any $T'\in(0,T)$ in the first inequality of~\eqref{hptecn}.
Moreover, we observe that \eqref{hpconvk} and \eqref{hptecn} imply that
$\ka\geq\kamin$, so that $\ka>0$ as a consequence.
Here is our first result.

\Bthm
\label{Convergenza}
Assume \Hpstruttura{\pier , {\gianni\Hpdati,}\ 
and let $(\thetae,\chie,\xie)$ be} the unique solution to problem~\Pbleps\
given by Theorem~\ref{Buonapositura}.
Moreover, assume \accorpa{hpconvk}{defka} and~\eqref{hptecn}.
Then, $(\thetae,\chie,\xie)$ converges in a proper topology
to the unique solution $(\theta,\chi,\xi)$ to problem~\Pbl\ satisfying \Regsoluz.
\Ethm

The topology mentioned in Theorem~\ref{Convergenza} 
will be clear from the proof we give in Section~\futuro\ref{Convergenza}
and is rather strong.
Provided that a much weaker topology is considered,
an error estimate can be proved.
We have indeed

\Bthm
\label{Errore}
Under the assumptions of Theorem~\ref{Convergenza}, the following estimate holds true
\Bsist
  && \norma{1*(\thetae-\theta)}_{\L\infty V}^2 + \norma{\chie-\chi}_{\L\infty H\cap\L2V}^2
  \non
  \\
  && \quad {}
  + \intQ (\ln\thetae-\ln\theta)(\thetae-\theta)
  + \intQ (\xie-\xi)(\chie-\chi)
  \non
  \\
  && \leq M \norma{1*\ke-\kalim}_{L^1(0,T)}
  \label{errore}
\Esist
for $\eps$ small enough,
where $M$ depends on the structure and the data, only.
\Ethm

\Brem
\label{Remnucleo}
We observe that a sufficient condition for {\pier\eqref{hptecn}$_2$} is that
$\ke$ has the form $\ke(t)=\eps^{-1}\hat k(t/\eps)$
with $\hat k\in L^1(0,+\infty)$.
Moreover, a sufficient condition for {\pier\eqref{hptecn}$_1$} to hold is that
$\ke$~is a positive type kernel, i.e.,
\Beq
  \ioT (\ke*v)(t) \, v(t) \,dt \geq 0
  \quad \hbox{for every $v\in L^2(0,T)$}.
  \label{tipopos}
\Eeq
In such a case, one can take $\kamin=\kaz$, indeed. {\elena The fact that the kernel is {\pier of positive type is actually in accordance with the assumptions required on $k$ to ensure thermodynamical consistency of
the model~\cite{BCF}}}.
We remark that a sufficient condition for a kernel $k$ to be of positive type
is the following (see, e.g., {\pier \cite[Prop.~4.1, p.~237]{Barbu} or \cite{GriLonSta}}):
$k$~is smooth, nonnegative, decreasing, and convex.
So, any positive multiple of $\exp(-t/\eps)$ plays the role
and the kernel given by
$\ke(t):=(\kalim/\eps)\exp(-t/\eps)$
is a prototype for both \eqref{hpconvk} and~\eqref{hptecn}
since $1*\ke$ converges to $\kalim$ strongly in~$L^1(0,T)$.
More generally,
assume that $\ke$ can be split as $\ke=p_\eps+r_\eps$,
where $p_\eps$ is of positive type and $r_\eps$ is a remainder.
If
\Bsist
  && 1*p_\eps \to \ka_1'
  \aand
  1*r_\eps \to \ka_2'
  \quad \hbox{strongly in $L^1(0,T)$}
  \non
  \\
  && \norma{r_\eps}_{L^1(0,T)} \leq {\pier\eta_0}
  \quad \hbox{for every $\eps>0$}
  \non
\Esist
for some constants $\ka_1'$, $\ka_2'$, and ${\pier\eta_0}<\kaz$,
then both assumptions \eqref{hpconvk} and \eqref{hptecn} are still fulfilled.
Indeed, we can take
$\kalim=\ka_1'+\ka_2'$, clearly, 
and $\kamin=\kaz-{\pier\eta_0}$, since 
\Beq
  \ioT \! (\kaz v + \ke*v ) v \, dt
  \geq \ioT \! (\kaz v + r_\eps*v ) v \, dt
  \geq (\kaz - {\pier\eta_0}) \ioT \! v^2 \, dt
  \quad \hbox{for every $v\in L^2(0,T)$}
  \non
\Eeq
thanks to the \holder\ and Young inequalities (see also \eqref{young}).
Note that $\ka_1'\geq0$ since $p_\eps$ is of positive type,
while $\ka_2'$ can be any real constant.
\Erem

\Brem
\label{Dirac}
Assumption \eqref{hpconvk} is a well-defined reinforcement
of the condition roughly mentioned in the Introduction as
$\ke\to\kalim\delta$, where $\delta$ is the Dirac mass at the origin.
Indeed, if we introduce the Heaviside function~$H$ on~$(-\infty,T)$, i.e.,
$H(t)=0$ for $t<0$ and $H(t)=1$ for $t\in(0,T)$,
and the trivial extension $\ketriv$ of~$\ke$,
\eqref{hpconvk} reeds
\Beq
  H * \ketriv \to \kalim H 
  \quad \hbox{strongly in $L^1(-\infty,T)$}
  \non
\Eeq
with an obvious new meaning of the convolution.
By differentiating and observing that
$(H*\ketriv)'=\delta*\ketriv=\ketriv$,
we deduce that
\Beq
  \ketriv \to \kalim\delta
  \quad \hbox{in the sense of distributions on $(-\infty,T)$}
  \non
\Eeq
where $\delta$ is the actually well-defined Dirac mass at $0$ in the open set $(-\infty,T)$.
\Erem

\Brem
\label{kzerotau}
{\pier By checking the proofs in the next sections, the reader 
will be able to realize that our results 
can be suitably extended to the case of 
coefficients $\kappa_{0\tau}$ possibly 
depending on $\tau $, with boundedness and 
convergence properties as $\tau \searrow 0$.}
\Erem

We recall that $\Omega$ is bounded and smooth.
So, throughout the paper,
we owe to some \wk\ embeddings of Sobolev type,
namely
$V\subset\Lx p$ for $p\in[1,6]$,
together with the related Sobolev inequality
\Beq
  \norma v_p \leq C \normaV v
  \quad \hbox{for every $v\in V$ and $1\leq p \leq 6$}
  \label{sobolev}
\Eeq
and $\Wx{1,p}\subset\Cx0$ for $p>3$, together with 
\Beq
  \norma v_\infty \leq C_p \norma v_{\Wx{1,p}}
  \quad \hbox{for every $v\in\Wx{1,p}$ and $p>3$}.
  \label{sobolevbis}
\Eeq
In \eqref{sobolev}, $C$ depends only on~$\Omega$, 
while $C_p$ in \eqref{sobolevbis} depends also on~$p$.
In particular, the continuous embedding 
$W\subset\Wx{1,6}\subset\Cx0$ holds.
Some of the previous embeddings are in fact compact.
This is the case for $V\subset\Lq$ and $W\subset\Cx0$.
We note that also the embeddings $W\subset V$, $V\subset H$,
$\Vz\subset H$, and $H\subset\Vp$ are compact.
Moreover, we often account for the \wk\ Poincar\'e inequalities
\Bsist
  && \normaV v \leq C \normaH{\nabla v}
  \quad \hbox{for every $v\in\Vz $}
  \label{poinczero}
  \\
  && \normaV v \leq C \Bigl( \normaH{\nabla v} + \bigl| \textstyle\iO v \bigr| \Bigr)
  \quad \hbox{for every $v\in V$}
  \label{poincmedia}
\Esist
where $C$ depends only on~$\Omega$.
Furthermore, we repeatedly make use of the notation
\Beq
  Q_t := \Omega \times (0,t)
  \quad \hbox{for $t\in[0,T]$}
  \label{defQt}
\Eeq
and of \wk\ inequalities, namely, the \holder\ inequality and the elementary Young inequality{\pier :}
\Beq
  ab \leq \delta a^2 + \frac 1{4\delta} \, b^2
  \quad \hbox{for every $a,b\geq 0$ and $\delta>0$}.
  \label{eleyoung}
\Eeq
As far as properties of the convolution are concerned,
we take advantage of the elementary formulas 
(which hold whenever they make sense)
\Beq
  a*b=a(0)(1*b)+a_t*1*b
  \aand
  (a*b)_t=a(0)b+a_t*b
  \label{comode}
\Eeq
and of the \wk\ Young theorem
\Beq
  \norma{u*v}_{\L rX}\leq \norma u_{L^p(0,T)} \norma v_{\L qX}
  \label{young}
\Eeq
where $X$ is a Banach space, $1\leq p,q,r\leq\infty$, and $1/r=(1/p)+(1/q)-1$
(cf., e.g., \cite{GriLonSta}).
Finally, again throughout the paper,
we use a small-case italic $c$ for different constants, that
may only depend 
on~$\Omega$, the final time~$T$, the shape of the nonlinearities $\lambda$, $\beta$, $\sigma$, 
and the properties of the data involved in the statements at hand; 
a~notation like~$c_\delta$ signals a constant that depends also on the parameter~$\delta$. 
The reader should keep in mind that the meaning of $c$ and $c_\delta$ might
change from line to line and even in the same chain of inequalities, 
whereas those constants we need to refer to are always denoted by 
capital letters, just like $C$ in~\eqref{sobolev}.


\section{Auxiliary material}
\label{Auxiliary}
\setcounter{equation}{0}

This section contains a very short summary on the properties
of the harmonic extension $\thetaH$ of the boundary datum~$\thetaG$
(see~\eqref{defthetaH})
and a preliminary result dealing with a \generaliz ed version of the limit problem \Pbl.
The properties listed in the following {\pier proposition} 
will be extensively used in the sequel.{\pier \relax

\Bprop
\label{StimethetaH}
Assumptions \eqref{regthetaG}--\eqref{hpthetaG} yield
\Beq
\label{piertheH}
  \thetaH \in  \H1H \cap \W{1,1}\Linfty \cap \C0V .
\Eeq
More precisely, owing to the theory of harmonic functions, in particular to
the maximum principle, we~have~that
\Bsist
  && \norma\thetaH_{\L2V}
    \leq C \norma\thetaG_{\L2{\HxG{1/2}}},
  \non
  \\
  && \norma{\thetaH}_{\L1H}
    \leq C \norma{\thetaG}_{\L1{\HxG{-1/2}}},
  \non
  \\
  && \norma{\dt\thetaH}_{\L2H}
    \leq C \norma{\dt\thetaG}_{\L2{\HxG{-1/2}}},
  \non
  \\
  && \thetamin \leq \thetaH \leq \thetamax \quad \aeQ,
  \non
  \\
  && \norma{\dt\thetaH}_{\L1\Linfty}
    = \norma{\dt\thetaG}_{\L1{L^\infty(\Gamma)}} 
  \non
\Esist
where $C$ is a constant depending on~$\Omega$, only.
\Eprop
}

Now, in order to {\pier help the reader},
we sketch the outline of the proof of Theorem~\ref{Convergenza}
we are going to develop in the next section.
By accounting for a number of a~priori estimates and using \wk\ compactness results,
we derive that the family of solutions $(\thetae,\chie,\xie)$
converges (for~a subsequence)
to a \generaliz ed solution to problem~\Pbl,
in which $\ln\theta$ is understood in a non standard sense.
Next, in order to conclude that such a solution actually 
is the solution given by Corollary~\ref{Buonaposituralimite},
{\pier we prove a preliminary} well-posedness result for \generaliz ed solutions
(Theorem~\futuro\ref{Buonaposituragen}).
Therefore, we first have to introduce the ingredients that are needed to explain such a notion of solution.

We define a \generaliz ed logarithm by following~\cite[Def.~4.2]{GR}.
However, we confine ourselves to consider the functions $\theta\in\L2V$
such that $\theta=\thetaG$ on the boundary, i.e., $\theta\in\thetaH+\L2\Vz$.
First, we introduce the function $\psi:\erre\to(-\infty,+\infty]$
by setting
\Beq
  \psi(r) = r(\log r - 1) \quad \hbox{if $r>0$},\quad
  \psi(0) = 0,
  \quad \hbox{and} \quad
  \psi(r) = +\infty \quad \hbox{if $r<0$}.
  \label{defpsi}
\Eeq
Then, for $\theta\in\thetaH+\L2\Vz$, 
we term $\Ln\theta$ the set of $\zeta\in\L2\Vp$ satisfying
\Beq
  \< \zeta , v - \theta > + \intQ \psi(\theta) \leq \intQ \psi(v)
  \quad \hbox{for every $v\in\theta+\L2\Vz$} 
  \label{defLn}
\Eeq
where $\<\cpto,\cpto>$ stands for the duality pairing 
between $\L2\Vp$ and $\L2\Vz$.
It {\pier can be checked} (see \cite[Thm.~4.7]{GR})) that
$\theta>0$ \aeQ\ whenever $\Ln\theta$ is not empty.
Moreover, even though $\Ln\theta$ might contain elements 
that are not functions (they are just Radon measure in such a case), 
its definition actually \generaliz es the usual logarithm.
Indeed, for $\theta\in\thetaH+\L2\Vz$
we have (see \cite[Rem.~4.3]{GR})
\Beq
  \ln\theta \in \Ln\theta
  \quad \hbox{whenever} \quad 
  \hbox{$\theta>0$ \aeQ\ and $\ln\theta\in\LQ2$}. 
  \label{lnLn}
\Eeq
Furthermore, the \generaliz ed logarithm is related to the theory of subdifferentials
as follows.
We define the function $\Psi:\L2\Vz\to(-\infty,+\infty]$ {\pier by}
\Beq
  \Psi(v) := \intQ \psi(v+\thetaH)
  \quad \hbox{for $v\in\L2\Vz$} 
\label{defPsi}
\Eeq
being understood that the integral is infinite if $\psi(v+\thetaH)\not\in\LQ1$.
Then, $\Psi$~turns out to be convex proper and \lsc\ on $\L2\Vz$, so that 
its (possibly multivalued) subdifferential $\partial\Psi:\L2\Vz\to\L2\Vp$
is well-defined.
Precisely, we have 
(see \cite[Rem.~4.4]{GR} for details)
\Beq
  \Ln\theta = \partial\Psi(\theta-\thetaH)
  \quad \hbox{{\pier for} every $\theta\in\thetaH+\L2\Vz$}.
  \label{cnsLn}
\Eeq
At this point, we can state the \generaliz ed version of problem~\Pbl\ as follows.
We look for a quadruplet $(\theta,\zeta,\chi,\xi)$ satisfying
\eqref{regtheta}, \accorpa{regchi}{regxi}, and
\Bsist
  && \zeta \in \H1\Vp 
  \aand
  \zeta \in \Ln\theta
  \label{regzeta}
  \\
  && \dt \bigl( \zeta(t) + \lambda(\chi(t)) \bigr)
  - \ka \Delta\theta(t) = f(t)
  \quad \hbox{in $\Vp$ \aat}
  \qquad
  \label{primagen}
  \\
  && \dt\chi - \Delta\chi + \xi + \sigma'(\chi)
  = \lambda'(\chi) \, \theta
  \aand
  \xi \in \beta(\chi) 
  \quad \aeQ
  \qquad
  \label{secondagen}
  \\
  && \zeta(0) = \ln\thetaz
  \aand
  \chi(0) = \chiz \,.
  \label{cauchygen}
\Esist
\Accorpa\Pblgen regzeta cauchygen

The following result holds

\Bthm
\label{Buonaposituragen}
Assume that both \Hpstruttura\ and {\gianni\Hpdati}\ hold
and that $\ka>0$.
Then, {\elena problem \Pblgen} has a unique solution $(\theta,\zeta,\chi,\xi)$
satisfying \eqref{regtheta} {\pier and \accorpa{regchi}{regxi}.} 
Moreover,  
{\pier we have that}
\Beq
  \ln\theta \in \L\infty H
  \aand
  \zeta = \ln\theta .
  \label{bellazeta}
\Eeq
\Ethm

\Bdim
We first prove uniqueness.
Our proof closely follows \cite[Sect.~5]{BCFG1}.
However, we repeat at least a part of the argument for the reader's convenience.
We first observe that, if $(\theta,\zeta,\chi,\xi)$ is a \generaliz ed solution, 
by integrating \eqref{primagen} in time we obtain
\Beq
  \zeta + \lambda(\chi) - \ka \Delta (1*\theta)
  = \ln\thetaz + \lambda(\chiz) + 1*f .
  \label{intprimagen}
\Eeq
Now, we pick two solutions $(\theta_i,\zeta_i,\chi_i,\xi_i)$, $i=1,2$,
and set for convenience
\Beq
  \theta := \theta_1 - \theta_2 \,, \quad
  \zeta := \zeta_1 - \zeta_2 \,, \quad
  \chi := \chi_1 - \chi_2 \,, \aand
  \xi := \xi_1 - \xi_2 \,.
  \label{defdiff}
\Eeq
Now, we write \eqref{intprimagen} for both solutions,
{\gianni take the difference, and test it
by $\theta\chart$ in the duality $\L2\Vp$-$\L2\Vz$,
where $t\in(0,T)$ is arbitrary and $\chart$ is the characteristic function of the interval $(0,t)$}.
We observe that all the terms but one are in fact integrals.
Namely, we~have
\Beq
  \< \zeta , \theta\chart >
  + \frac \ka 2 \iO |\nabla (1*\theta)(t)|^2
  = - \intQt \bigl( \lambda(\chi_1) - \lambda(\chi_2) \bigr) \, \theta .
  \label{diffintprimagen}
\Eeq
At the same time, we write \eqref{secondagen} for both solutions,
take the difference, multiply the {\pier resulting} equality by~$\chi$,
and integrate over~$Q_t$.
We~{\pier obtain}
\Bsist
  && \frac 12 \iO |\chi(t)|^2
  + \intQt |\nabla\chi|^2
  + \intQt \xi\chi
  \non
  \\
  && = \intQt \bigl( \lambda'(\chi_1) \, \theta_1 - \lambda'(\chi_2) \, \theta_2 \bigr) \chi
  - \intQt {\pier \bigl( \sigma'(\chi_1) - \sigma'(\chi_2) \bigr) \chi} .
  \label{diffsecondagen}
\Esist
Finally, we add \eqref{diffsecondagen} to \eqref{diffintprimagen}
and {\pier proceed} exactly as in \cite{BCFG1}.
{\pier However, let us point out that, in view of Taylor's expansion
and the H\" older and Sobolev inequalities, we~have
\Bsist
  && \displaystyle \intQt
      \bigl\{ - 
        \bigl( \lambda(\chi_1) - \lambda(\chi_2) \bigr) \theta
        + \bigl(
          \lambda'(\chi_1) \theta_1 - \lambda'(\chi_2) \theta_2
          \bigr) \chi
      \bigr\}
  \cr
  && \displaystyle  =  \intQt \theta_1
      \graffe{
        \lambda(\chi_2) - \lambda(\chi_1) - \lambda'(\chi_1)( \chi_2 -\chi_1)
      } \cr
  && \displaystyle \qquad + \intQt \theta_2
      \graffe{
        \lambda(\chi_1) - \lambda(\chi_2) - \lambda'(\chi_2)( \chi_2 -\chi_1)
      }
  \cr
  && \displaystyle  \leq c \intQt (\theta_1 + \theta_2) |\chi|^2
    \leq c \int_0^t \norma{\theta_1(s) + \theta_2(s)}_4
      \norma{\chi(s)}_4 \norma{\chi(s)}_2 \, ds
  \cr
  && \displaystyle \leq c \int_0^t \normaV{\theta_1(s) + \theta_2(s)}
      \bigl( \normaH{\nabla\chi(s)} + \normaH{\chi(s)} \bigr)
      \normaH{\chi(s)} \, ds
  \cr
  && \displaystyle  \leq \frac14 \intQt |\nabla\chi|^2
    + c \int_0^t
        {\gianni\tonde{ \normaV{\theta_1(s)}^2 + \normaV{\theta_2(s)}^2 }}
        \normaH{\chi(s)}^2 \, ds .
 \nonumber
\Esist

Therefore,} we derive~that
\Bsist
  && \< \zeta , \theta\chart >
  + \frac \ka 2 \iO |\nabla (1*\theta)(t)|^2
  + \frac 12 \iO |\chi(t)|^2
  + \intQt |\nabla\chi|^2
  + \intQt \xi\chi
  \non
  \\
  && {}\leq {}{\pier 
     c \iot
       \tonde{1+ \normaV{\theta_1(s)}^2 + \normaV{\theta_2(s)}^2 }
        \normaH{\chi(s)}^2 \, ds 
        + \frac 14 \intQt |\nabla\chi|^2.}
  \non
\Esist
{\pier As the last term can be easily controlled by the 
\lhs,}
what remains to observe is that all the terms on the \lhs\ are non negative.
The integral containing $\xi$ is non negative by monotonicity.
Let us deal with the duality term.
For $i=1,2$ we have $\zeta_i\in\Ln\theta_i$.
Due to~\eqref{defLn}, this means that
\Beq
  \< \zeta_i , v_i - \theta_i > + \intQ \psi(\theta_i) \leq \intQ \psi(v_i)
  \quad \hbox{for every $v_i\in{\pier \theta_i}+\L2\Vz$ and $i=1,2$} .
  \non
\Eeq
Now, we choose the admissible functions 
$v_1=\theta_1-\theta\chart$ and $v_2=\theta_2+\theta\chart$,
we sum up and split the integrals.
We have
\Bsist
  && \< \zeta , - \theta \chart >
  + \intQt \bigl( \psi(\theta_1) + \psi(\theta_2) \bigr)
  + \int_{Q\meno Q_t} \bigl( \psi(\theta_1) + \psi(\theta_2) \bigr)
  \non
  \\
  && \leq \intQt \bigl( \psi(v_1) + \psi(v_2) \bigr)
  + \int_{Q\meno Q_t} \bigl( \psi(v_1) + \psi(v_2) \bigr).
  \non
\Esist
As $v_1=\theta_2$ in $Q_t$ and $v_1=\theta_1$ in $Q\meno Q_t$
and similarly for~$v_2$, all the integrals are finite and cancel out.
We deduce that
$\<\zeta,-\theta\chart>\leq0$, i.e., what we wanted to prove.
At this point, we can apply the Gronwall lemma (see, e.g., \cite[pp.156-157]{Brezis})
and obtain, in particular, that $\chi=0$ and $\nabla(1*\theta)=0$ \aeQ.
As $1*\theta$ is $\Vz$-valued, this implies that $1*\theta=0$ \aeQ,
whence also $\theta=0$ \aeQ.
All this means that $\theta_1=\theta_2$ and $\chi_1=\chi_2$.
By comparison in \eqref{intprimagen} and~\eqref{secondagen},
we conclude that $\zeta_1=\zeta_2$ and $\xi_1=\xi_2$ as well.

Once uniqueness of the \generaliz ed solution is proved,
we can easily conclude.
Indeed, our assumptions allow us to apply Corollary~\ref{Buonaposituralimite}.
Hence, a~solution exists in the strong sense, i.e., satisfying the regularity requirements \Regsoluz.
On the other hand, such a solution is also a \generaliz ed solution due to~\eqref{lnLn}.
Finally, it satisfies~\eqref{bellazeta}.
\Edim


\section{Proofs of Theorems \ref{Convergenza} and \ref{Errore}}
\label{Proofs}
\setcounter{equation}{0}

The argument we follow for {\pier our first proof} 
uses compactness and monotonicity methods.
So, we start estimating.
However, we often proceed formally for the sake of simplicity.
The correct procedure could be based on performing similar estimates 
on the solution of {\pier some approximating problem. One approximation is} constructed in~\cite{BCFG1} 
{\pier and} depends on the parameter~$\epsilon${\pier : 
the solution} is smoother than 
the solution to {\pier the problem 
we are dealing with}
({\elena actually} the~limit as $\epsilon\todx0$ {\pier keeps} such estimates).
Furthermore, in order to simplify the notation,
we often avoid the subscript $\eps$ {\elena (on the solutions)} during the calculation
and restore it just at the end of each estimate.

\step First a priori estimate

We would like testing \eqref{primaeps} by 
\Beq
  \ue := \thetae - \thetaH
  \label{defueps}
\Eeq
in the duality $\Vp$-$\Vz$ and integrate over $(0,t)$,
where $t\in(0,T)$ is arbitrary.
However, we proceed formally, as just said.
In particular, we behave as if the logarithmic term were smoother.
At the same time we multiply \eqref{secondaeps} by $\dt\chie$
and integrate over~$Q_t$ (see~\eqref{defQt}).
Finally, we sum up and remark that the terms containing $\dt\chie$ partially cancel.
Hence, by~avoiding some subscripts in the notation for a while
and adding the same integral $\iO|\chi(t)|^2$ to both sides for convenience{\pier ,
we obtain}
\Bsist
  && \iO \theta(t)
  - \intQt \dt(\ln\theta) \, \thetaH
  + \intQt \Bigl( {\pier \kaz |\nabla u|^2 + \bigl( \ke*\nabla u \bigr) \cdot 
     \nabla u } \Bigr)
  \non
  \\
  && \quad {}
  + \intQt |\dt\chi|^2 
  + \frac 12 \iO |\nabla\chi(t)|^2
  + \iO \Beta(\chi(t))
  + \iO |\chi(t)|^2
  \non
  \\
  && = \intQt \dt\lambda(\chi) \, \thetaH
  {\pier {}- \intQt \Bigl( \kaz \nabla \thetaH + \bigl( \ke*\nabla \thetaH \bigr)
    \Bigr) \cdot \nabla u
  + \intQt f u}  
  \non
  \\
  && \quad {}+ \iO \bigl( |\chi(t)|^2 - \sigma(\chi(t)) \bigr) 
  + \iO \thetaz
  + \frac 12 \iO |\nabla\chiz|^2
  + \iO \bigl( \Beta+\sigma \bigr) (\chiz) .
  \qquad
  \label{perprima}
\Esist
{\elena Note that we have $\theta>0$ and {\pier can}  make use of the chain rule for subdifferentials}. 
Now, we recall that $\Beta$ is nonnegative (cf.~\eqref{hpBeta})
and treat each of the non-trivial terms, separately.
We integrate the second integral on the \lhs\ by parts 
with respect to time and~get
\Bsist
  && \intQt \dt(\ln\theta) \, \thetaH
  = \iO {\pier \thetaH (t)}\, \ln\theta(t) - 
   \iO {\pier \thetaH (0)}\,\ln\thetaz
  - \intQt \ln\theta \, \dt\thetaH 
  \non
  \\
  && \leq \iO {\pier \thetaH (t)}\,\ln^+\theta(t) 
  - \iO {\pier \thetaH (t)}\,\ln^-\theta(t) 
  + \intQt |\ln\theta| \, |\dt\thetaH| + c .
  \non
\Esist
Hence, by {\pier recalling Proposition~\ref{StimethetaH} and observing that 
$r-\theta^* \ln^+r\geq (r/2)-c $}\ for every $r>0$, 
we deduce~that
\Bsist
  && \iO \theta(t)
  - \intQt \dt(\ln\theta) \, \thetaH
  \non
  \\
  &&{\pier \geq } \iO \Bigl( \theta(t) - {\pier \theta^*}\ln^+\theta(t) 
  + {\pier \theta_*} \ln^-\theta(t) \Bigr)
  - \intQt {\gianni |\ln\theta| \, |\dt\thetaH|} - c
  \non
  \\
  && \geq \iO \Bigl( \frac 12 \, \theta(t) + {\pier \theta_*} \ln^-\theta(t) \Bigr)
  - \intQt |\ln\theta| \, |\dt\thetaH| - c .
  \non
\Esist
On the other hand, we notice that $|\ln r|\leq c((r/2)+{\pier \theta_*}\ln ^-r)$ for $r>0$, 
so~that
\Beq
  \intQt |\ln\theta| \, |\dt\thetaH|
  \leq c \intQt \bigl( (\theta/2) + {\pier \theta_*}\ln^-\theta \bigr) \, |\dt\thetaH|
  \non
\Eeq
and the last integral can be treated on the \rhs\ via the Gronwall lemma
since $\dt\thetaH\in\L1\Linfty$.
Next, {\pier thanks} to \eqref{hptecn} {\gianni and to~\eqref{poinczero}, we infer~that}
\Beq
  \intQt \Bigl( \kaz |\nabla{\pier u} |^2 + \bigl( \ke*\nabla{\pier u} 
  \bigr) \cdot \nabla{\pier u} \Bigr)
  \geq \kamin \intQt |\nabla{\pier u}|^2 
  \geq {\gianni \frac \kamin C} \int_0^t \| u (s)\|_V^2 \ds 
  \non
\Eeq
{\gianni where $C$ is the constant in~\eqref{poinczero}}.
Now, let us deal with the \rhs. 
{\pier With the help of \eqref{hptecn}, {\gianni the Young theorem~\eqref{young}, and \eqref{hpf}}, we} immediately have
\Beq
- \intQt \Bigl( \kaz \nabla \thetaH + \bigl( \ke*\nabla \thetaH \bigr)
    \Bigr) \cdot \nabla u
  + \intQt f u \leq {\pier {}c + {\gianni c} \int_0^t \| u (s)\|_V^2 \ds} .
  \non
\Eeq
Moreover, we {\pier observe that}
\Beq
  \intQt \dt\lambda(\chi) \, \thetaH
  = \intQt \lambda'(\chi) \, \dt\chi \, \thetaH
  \leq \frac 14 \intQt |\dt\chi|^2 + c \intQt |\chi|^2 + c
  \non
\Eeq
since $|\lambda'(r)|\leq c(1+|r|)$ by \eqref{hpls} and $0\leq\thetaH\leq\thetamax$ by {\pier Proposition}~\ref{StimethetaH}.
Finally, we observe that \eqref{hpls} also yields $|\sigma(r)|\leq c(1+r^2)$ for every $r$
and deduce~that
\Bsist
  && \iO \bigl( |\chi(t)|^2 - \sigma(\chi(t)) \bigr)
  \leq c + c \iO |\chi(t)|^2
  \leq c + c {\gianni \iO \Bigl( |\chiz|^2 + \iot 2 \chi(s) \, \dt\chi(s) \, ds \Bigr)}
  \non
  \\
  && \leq {\pier {}c{}+{} }\frac 14 \intQt |\dt\chi|^2 + c \intQt |\chi|^2 .
  \non
\Esist
By combining all the estimates we have derived with~\eqref{perprima},
applying the Gronwall lemma, and owing to the Poincar\'e inequality~\eqref{poincmedia},
we conclude~that
\Beq
  \norma\thetae_{\L\infty\Luno\cap\L2 V}
  + \norma{\ln\thetae}_{\L\infty\Luno}
  + \norma\chie_{{\pier\H1H\cap\L\infty V}}
  \leq c 
  \label{primastima}
\Eeq
besides an estimate for $\Beta(\chie)$ in $\L\infty\Luno$.

\step Second a priori estimate

We write \eqref{secondaeps} in the form of a nonlinear monotone elliptic equation,
namely
\Beq
  -\Delta\chie + \xie = \lambda'(\chie) \, \thetae - \dt\chie - \sigma'(\chie)
  \aand 
  \xie \in \beta(\chie)
  \label{perstimaxi}
\Eeq
and notice that each term on the \rhs\ of \eqref{perstimaxi}
is bounded in $\L2H$ by \eqref{hpls} and~\eqref{primastima}.
{\pier Concerning the first term,}
notice that $\chie$ and $\thetae$ are bounded 
in $\L\infty{\Lx4}$ and $\L2{\Lx4}$, respectively, 
due to the Sobolev inequality~\eqref{sobolev}.
Then, a~quite standard argument 
(formally test \eqref{perstimaxi} by either $-\Delta\chie$ or~$\xie$
in order to estimate both of them
and then use the regularity theory for elliptic equations)
yields
\Beq
  \norma\chi_{\L2W} + \norma\xie_{\L2H}
  \leq c .
  \label{stimaxi} 
\Eeq

\step Third a priori estimate

We want to estimate $\dt\thetae$ in $\L1{\Wx{-1,q}}$ for some $q>1$
satisfying $\Ldue\subset\Wx{-1,q}$,
and the choice $q=4/3$ will work.
Therefore, by proceeding formally, 
we take any $v\in\Wxz{1,4}$ satisfying $\norma v_{\Wx{1,4}}\leq1$
and test \eqref{primaeps} written at {\pier almost} any time~$t$ by $\thetae(t)\,v$, {\elena which is a good test function 
belonging
to $V_0$, due to \eqref{sobolev}-\eqref{sobolevbis} and \eqref{primastima}}.
We obtain {\elena (the first integral is intended as a duality pairing)}
\Bsist
  & \displaystyle \iO \dt\thetae(t) \, v
  & = \iO \bigl\{ f(t) - \lambda'(\chie(t)) \dt\chie(t) \bigr\} \thetae(t) \, v
  \non
  \\
  && - \iO \kaz \nabla\thetae(t) \cdot \nabla(\thetae(t) \, v)
  - \iO (\ke*\nabla\thetae)(t) \cdot \nabla(\thetae(t) \, v) .
  \label{perterza}
\Esist
We simplify the notation by dropping the time~$t$ (and the subscript~$\eps$, as usual) for a while
and estimate each term on the \rhs\ of~\eqref{perterza}, separately.
We account for the Lipschitz continuity of~$\lambda'$ (see~\eqref{hpls})
and for the \holder\ and Sobolev inequality~\eqref{sobolev}.
Moreover, we observe that $\norma v_\infty\leq c$ thanks to~\eqref{sobolevbis}. 
We~have 
\Bsist
  && \iO \bigl\{ f - \lambda'(\chi) \dt\chi \bigr\} \theta v
  \leq \norma f_2 \norma\theta_2 \norma v_\infty
  + \norma{\lambda'(\chi)}_4 \norma{\dt\chi}_2 \norma\theta_4 \norma v_\infty
  \non
  \\
  && \qquad {}
  \leq c \norma f_2 \norma\theta_2 
  + c (1 + \norma\chi_4) \norma{\dt\chi}_2 \norma\theta_4 \,.
  \non
  \\
  && - \iO \kaz \nabla\theta \cdot \nabla(\theta \, v)
  \leq {\gianni \kaz \,} \norma{\nabla\theta}_2^2 \norma v_\infty
  + \kaz \norma{\nabla\theta}_2 \norma\theta_4 \norma{\nabla v}_4 
  \non
  \\
  && \qquad {}
  \leq c \norma{\nabla\theta}_2^2
  + c \norma{\nabla\theta}_2 \norma\theta_4 \,.
  \non
  \\
  && - \iO (\ke*\nabla\thetae) \cdot \nabla(\thetae \, v) 
  \leq \norma{\ke*\nabla\theta}_2 \bigl( \norma{\nabla\theta}_2 \norma v_\infty + \norma\theta_4 \norma{\nabla v}_4 \bigr)
  \non
  \\
  && \qquad {}
  \leq c \norma{\ke*\nabla\theta}_2 \bigl( \norma{\nabla\theta}_2 + \norma\theta_4 \bigr) .
  \non
\Esist
We term $C$ the maximum of the values of the above constant~$c$'s,
for clarity.
Then, we first collect the estimates just obtained and \eqref{perterza}.
Finally, we take the supremum with respect to~$v\in\Wxz{1,4}$ 
under the {\elena constraint}  $\norma v_{\Wx{1,4}}\leq 1$. 
We conclude~that
\Beq
  \norma{\dt\thetae(t)}_{\Wx{-1,4/3}}
  \leq C \, \phie(t)
  \quad \aat
  \non
\Eeq
where $\phie:(0,T)\to\erre$ is defined by 
\Bsist
  && \phie(t) :=
  \norma{f(t)}_2 \norma{\thetae(t)}_2 
  + (1 + \norma{\chie(t)}_4) \norma{\dt\chie(t)}_2 \norma{\thetae(t)}_4 
  + \norma{\nabla\thetae(t)}_2^2
  \non
  \\
  && \qquad  {}
  + \norma{\nabla\thetae(t)}_2 \norma{\thetae(t)}_4 
  + \norma{(\ke*\nabla\thetae)(t)}_2 \bigl( \norma{\nabla\thetae(t)}_2 + \norma{\thetae(t)}_4 \bigr) .
  \non
\Esist
Therefore, the estimate
\Beq
  \norma{\dt\thetae}_{\L1{\Wx{-1,4/3}}} \leq c 
  \label{terzastima}
\Eeq
follows once we prove that $\phie$ is bounded in~$L^1(0,T)$.
In view of the previous estimates \eqref{primastima}, \eqref{stimaxi},
and of the Sobolev inequality,
we see that the only trouble could come from the term containing the convolution.
By owing to the Young theorem (see~\eqref{young})
and to the Sobolev inequality once more,
we have
\Bsist
  && \ioT \norma{(\ke*\nabla\thetae)(t)}_2 \bigl( \norma{\nabla\thetae(t)}_2 + \norma{\thetae(t)}_4 \bigr) \, dt
  \non
  \\
  && \leq c \norma{\ke*\nabla\thetae}_{\L2H} \norma\thetae_{\L2V}
  \leq c \norma\ke_{L^1(0,T)} \norma\thetae_{\L2V}^2 \,. 
  \non
\Esist
Now, we recall {\pier \eqref{hptecn}$_2$} and \eqref{primastima}
and conclude that $\phie$ is bounded in $L^1(0,T)$.
Therefore, \eqref{terzastima} is established.

\step Fourth a priori estimate

By testing \eqref{primaeps} by any $v\in\L2\Vz$
and integrating over~$(0,T)$,
we deduce~that
\Beq
  \Bigl| \ioT \< \dt\ln\thetae(t) \, , v(t) > \, dt \Bigr|
  \leq {\pier c \,} M_\eps  \norma v_{\L2V}
  \quad \hbox{for every $v\in\L2\Vz$}
  \non
\Eeq
where we have set
\Bsist
  M_\eps
  := \kaz \norma{\nabla\thetae}_{\L2H} 
  + \norma{\ke*\nabla\thetae}_{\L2H}
  + \norma f_{\L2H} \non \\
  {}+ {\pier \tonde{1+ \norma{\chie}_{\L\infty{L^4(\Omega)}} }} \,  \norma{{\gianni \dt\chie}}_{\L2H} \,.
  \label{perstimadtln}
\Esist
Thus, the estimate
\Beq
  \norma{\dt\ln\thetae}_{{\elena \L2{\Vz^*}}} \leq c
  \label{stimadtln}
\Eeq
follows whenever we prove that $M_\eps\leq c$.
So, let us examine each term of \eqref{perstimadtln}
but the third one, of course, by accounting for~\eqref{primastima}.
The first and last ones are bounded by \eqref{hptecn} and~\eqref{hpls}.
For the second term, we use the Young inequality \eqref{young} and \eqref{hptecn}
and obtain
\Beq
  \norma{\ke*\nabla\thetae}_{\L2H}
  \leq \norma\ke_{L^1(0,T)} \norma{\nabla\thetae}_{\L2H}
  \leq c.
  \non
\Eeq
Hence, \eqref{stimadtln} is established.

\step Convergence and conclusion

From estimates \eqref{primastima} and \eqref{stimadtln}
we derive the following {\gianni convergence}
\Bsist
  & \thetae \to \theta
  & \quad
  {\pier \hbox{weakly in $\L2V $}}
  \label{convtheta}
  \\
  & \ue  \to u 
  & \quad
  \hbox{weakly in $\L2\Vz$}
  \label{convu}
  \\
  & \chie \to \chi
  & \quad 
  \hbox{weakly star in $\H1H\cap\L\infty V\cap\L2W$}
  \quad
  \label{convchi}
  \\
  & \xie \to \xi
  & \quad
  \hbox{weakly in $\L2H$}
  \label{convxi}
  \\
  & \ln\thetae \to \zeta
  & \quad
  \hbox{{\pier weakly in} {\elena $\H1{\Vz^*}$}}
  \label{convln}
\Esist
for suitable functions $\theta,\chi,\xi,\zeta$,
and $u:=\theta-\thetaH$,
possibly for a subsequence $\eps=\eps_n\todx0$.
We are going to {\pier show} that the triplet $(\theta,\chi,\xi)$
is a solution to problem \Pbl.
Once this is proved, the whole family $(\thetae,\chie,\xie)$ converges
(in~the above topology) to~$(\theta,\chi,\xi)$, since the solution to problem \Pbl\ is unique.
Now, we recall Theorem~\ref{Buonaposituragen} and observe that
it is sufficient to show that the quadruplet
$(\theta,\zeta,\chi,\xi)$ is a \generaliz ed solution,
i.e., it solves \Pblgen.
So, we just prove this fact.
The regularity requirements \eqref{regtheta}, \accorpa{regchi}{regxi}, 
and~\eqref{regzeta} are already clear.
First, we observe that \eqref{convchi} and \eqref{convln}
imply at least weak convergence in $\C0H$ and $\C0\Vp$, respectively,
whence $\zeta$ and $\chi$ satisfy the Cauchy conditions
$\zeta(0)=\ln\thetaz$ and $\chi(0)=\chiz$.
Furthermore, as the embeddings $W\subset V$, $V\subset\Lq$, 
and {\pier $V\subset H$} {\gianni are} compact,
{\pier we can apply, e.g., 
\cite[Sect.~8, Cor.~4]{Simon} 
and deduce that \eqref{convchi}, \eqref{convtheta}}  
and \eqref{terzastima} imply
\Bsist
  && \chie \to \chi
  \quad \hbox{strongly in $\L2V\cap\C0\Lq$}
  \label{convfortechi}
  \\
  && \thetae \to \theta
  \quad \hbox{strongly in $\L2H$}.
  \label{convfortetheta}
\Esist
Moreover, we can even assume that 
\Beq
  \chie\to\chi 
  \aand
  \thetae\to\theta
  \quad \aeQ .
  \label{convqo}
\Eeq
In particular, by noting that \eqref{hpls} implies
Lipschitz continuity for $\lambda'$ and $\sigma'$ 
and the estimate (via Taylor's formula)
\Beq
  |\lambda(\chie) - \lambda(\chi)|
  \leq c (1+|\chi|) \, |\chie-\chi|
  + c |\chie-\chi|^2 
  \non
\Eeq
we infer that
\Bsist
  & \lambda'(\chie) \to \lambda'(\chi)
  & \quad \hbox{strongly in~$\C0{\Lx4}$ and \aeQ}
  \label{convlap}
  \\
  & \sigma'(\chie) \to \sigma'(\chi)
  & \quad \hbox{strongly in~$\C0{\Lx4}$ and \aeQ}
  \label{convsip}
  \\
  & \lambda(\chie) \to \lambda(\chi)
  & \quad \hbox{strongly in $\C0H$ and \aeQ} .
  \label{convla}
\Esist
Next, we observe that \eqref{convtheta} implies that
$\Delta\thetae$ converges to $\Delta\theta$ weakly in $\L2\Vp$.
By~\accorpa{hpconvk}{defka} and the Young theorem (see~\eqref{young}), we infer that
\Beq
  (\kaz + 1 * \ke) * \Delta\thetae
  \to \ka * \Delta\theta
  \quad \hbox{weakly in $\L2\Vp$}
  \non
\Eeq
and we can take the limit in the integrated version of~\eqref{primaeps}.
Namely, we obtain
\Beq
  \zeta + \lambda(\chi) - \ka * \Delta\theta = 1*f + \ln\thetaz + \lambda(\chiz) 
  \label{perdopo}
\Eeq
and we conclude that the limit functions we have constructed satisfy {\pier \eqref{primagen} and}
\Beq
 {\pier  \dt\chi - \Delta\chi + \xi + \sigma'(\chi) = \lambda'(\chi) \, \theta \quad \aeQ.}
  \non
\Eeq
Thus, it just remains to prove that
$\xi\in\beta(\chi)$ \aeQ\ and $\zeta\in\Ln\theta$.
The first claim immediately follows by applying, e.g., \cite[Lemma~1.3, p.~42]{Barbu}
in the framework of maximal monotone operators in~$\LQ2$
by accounting for \eqref{convxi} and~\eqref{convfortechi}.
On the contrary, the second claim needs much more work.
We use the framework of the maximal monotone graphs in $\L2\Vz\times\L2\Vp$
and consider the subdifferential $\partial\Psi$ of the function~$\Psi$ given by~\eqref{defPsi},
which is related to the multivalued operator $\Ln$ by~\eqref{cnsLn}.
So, we have to prove that $\zeta\in\partial\Psi(u)$
by starting from $\ln\thetae\in\partial\Psi(\ue)$ for $\eps>0$
(see~\eqref{lnLn}) and the weak {\gianni convergence} given by \eqref{convu} and~\eqref{convln}.
It is well known (see  once more, e.g., \cite[Lemma~1.3, p.~42]{Barbu})
that a condition that allows to conclude is the following
\Beq
  \limsup_{\eps\todx0} \, \< \ln\thetae , \ue >
  \leq \< \zeta , u > , 
  \quad \hbox{that is,} \quad
  \limsup_{\eps\todx0} \, \< \ln\thetae - \zeta , \thetae - \theta >
  \leq 0.
  \label{perconcludere}
\Eeq
In order to prove \eqref{perconcludere},
we consider the integrated version of~\eqref{secondaeps} and equation~\eqref{perdopo}.
We take the difference and test the equality obtained by 
${\pier \thetae-\theta}$.
We have
\Bsist
  && \< \ln\thetae - \zeta , \thetae - \theta >
  = - \kaz \intQ \nabla \bigl( 1*(\thetae-\theta) \bigr) \cdot \nabla (\thetae-\theta)
  \non
  \\
  && \quad {}
  - \intQ \bigl( 1*\ke*\nabla\thetae - \kalim * \nabla\theta \bigr) \cdot \nabla(\thetae-\theta)
  - \intQ \bigl( \lambda(\chie) - \lambda(\chi) \bigr) (\thetae-\theta) 
  \non
  \\
  && =  - (\kaz+\kalim) {\gianni\intQ} \nabla \bigl( 1*(\thetae-\theta) \bigr) \cdot \nabla (\thetae-\theta)
  \non
  \\
  && \quad {}
  - \intQ \bigl( (1*\ke - \kalim) *\nabla\thetae \bigr) \cdot \nabla(\thetae-\theta)
  - \intQ \bigl( \lambda(\chie) - \lambda(\chi) \bigr) (\thetae-\theta) .
  \non
\Esist
Now, the first integral of the last chain is nonnegative and $\kaz+\kalim=\ka>0$.
Next, the last integral tends to zero by~\eqref{convfortetheta} and~\eqref{convla}.
Finally, the middle term is estimated by the \holder\ inequality, the Young theorem,
and \eqref{primastima} this way
\Bsist
  && - \intQ \bigl( (1*\ke - \kalim) *\nabla\thetae \bigr) \cdot \nabla(\thetae-\theta)
  \non
  \\
  && \leq \norma{(1*\ke - \kalim) *\nabla\thetae}_{\L2H} \norma{\nabla(\thetae-\theta)}_{\L2H}
  \non
  \\
  && \leq \norma{1*\ke - \kalim}_{L^1(0,T)} \norma{\nabla\thetae}_{\L2H} \norma{\nabla(\thetae-\theta)}_{\L2H}
  \non
  \\
  && \leq c \norma{1*\ke - \kalim}_{L^1(0,T)} .
  \label{perdimerrore}
\Esist
Hence, by recalling {\pier \eqref{hpconvk}},
we conclude that \eqref{perconcludere} holds.

\step
Proof of Theorem \ref{Errore}

By arguing as in the last part of the previous proof,
we consider the integrated version of the equations for temperature
and test the difference by $\thetae-\theta$.
However, in the present situation we already know that $\zeta=\ln\theta$
and can integrate over $Q_t$ rather than~$Q$.
Thus, a~quite similar calculation yields
\Bsist
  && \intQt (\ln\thetae - \ln\theta) (\thetae-\theta)
  + \ka \intQt \nabla \bigl( 1*(\thetae-\theta) \bigr) \cdot \nabla (\thetae-\theta)
  \non
  \\
  && = - \intQt \bigl( (1*\ke - \kalim) *\nabla\thetae \bigr) \cdot \nabla(\thetae-\theta)
  - \intQt \bigl( \lambda(\chie) - \lambda(\chi) \bigr) (\thetae-\theta) .
  \qquad
  \label{pererroreA}
\Esist
At the same time, we multiply the difference between \eqref{secondaeps} and \eqref{seconda}
by $\chie-\chi$ and integrate over~$Q_t$.
We obtain
\Bsist
  && \frac 12 \iO |(\chie-\chi)(t)|^2
  + \intQt |\nabla(\chie-\chi)|^2
  + \intQt (\xie-\xi)(\chie-\chi)
  \non
  \\
  && = \intQt \bigl( \lambda'(\chie) \, \thetae - \lambda'(\chi) \, \theta \bigr) (\chie-\chi)
  - \intQt \bigl( \sigma'(\chie) - \sigma'(\chi) \bigr) (\chie-\chi) .
  \label{pererroreB}
\Esist
At this point, we sum \eqref{pererroreB} to~\eqref{pererroreA},
and it is clear that all the terms on the \lhs\ are nonnegative.
Thus, we estimate each term on the \rhs.
As far as the term containing the convolution kernel is concerned,
we can repeat the argument that led to~\eqref{perdimerrore}.
Thus, we~have
\Beq
  - \intQ \bigl( (1*\ke - \kalim) *\nabla\thetae \bigr) \cdot \nabla(\thetae-\theta)
  \leq c \norma{1*\ke - \kalim}_{L^1(0,T)} .
  \non
\Eeq
The sum of all the terms involving $\lambda$ and $\lambda'$
can be first transformed and then estimated as follows
(we~use the Taylor formula and~\eqref{hpls}, besides standard inequalities, as usual)
\Bsist
  && - \intQt \thetae
      \graffe{
        \lambda(\chie) - \lambda(\chi) - \lambda'(\chie) (\chie-\chi)
      }
    - \intQt \theta
      \graffe{
        \lambda(\chi) - \lambda(\chie) + \lambda'(\chi) (\chie-\chi)
      }
  \non
  \\
  && \leq c \intQt (\thetae + \theta) \, |\chie-\chi|^2
  \non
  \\
  && \leq c \iot \norma{\thetae(s)+\theta(s)}_4 \norma{\chie(s)-\chi(s)}_4 \norma{\chie(s)-\chi(s)}_2 \, ds
  \non
  \\
  && \leq c \iot \normaV{\thetae(s)+\theta(s)} 
       \bigl( \normaH{\chie(s)-\chi(s)} + \normaH{\nabla(\chie(s)-\chi(s))} \bigr) \normaH{\chie(s)-\chi(s)} \, ds
  \non
  \\
  && \leq \frac 12 \intQt |\nabla(\chie-\chi)|^2
  + c \iot  \bigl( 1 + \normaV{\thetae(s)+\theta(s)}^2 \bigr) \normaH{\chie(s)-\chi(s)}^2 \, ds .
  \non
\Esist
As the integral involving $\sigma'$ can be treated in a trivial way due to~\eqref{hpls},
we can apply the Gronwall lemma (see, e.g., {\pier \cite[p.~156]{Brezis}})
and infer that the \lhs\ of \eqref{errore} is bounded~by
\Beq
  c \norma{1*\ke - \kalim}_{L^1(0,T)} \, \exp \Bigl( c \ioT \bigl( 1 + \normaV{\thetae(t)+\theta(t)}^2 \bigr) \, dt \Bigr).
  \non
\Eeq
As the last integral is bounded by a constant thanks to~\eqref{primastima},
inequality \eqref{errore} follows.



\vspace{3truemm}

\Begin{thebibliography}{10}

{\elenabis \bibitem{librofabrizio}
G. Amendola, M. Fabrizio, J.M. Golden, 
``Thermodynamics of materials with memory. Theory and applications'', 
Springer, New York, 2012.}

\bibitem{Barbu}
V. Barbu,
``Nonlinear semigroups and differential equations in Banach spaces'',
Noordhoff International Publishing, Leyden, 1976.

\bibitem{BCF}
E. Bonetti, P. Colli, M. Fr\'emond,
A phase field model with thermal memory governed by the entropy balance,
{\it Math. Models Methods Appl. Sci.}
{\bf 13} (2003) 1565-1588.

\bibitem{BCFG2}
E. Bonetti, P. Colli, M. Fabrizio, G. Gilardi,
Modelling and long-time behaviour for phase transitions 
with entropy balance and thermal memory conductivity,
{\it Discrete Contin. Dyn. Syst. Ser.~B} {\bf 6} (2006) 1001-1026. 

\bibitem{BCFG1}
E. Bonetti, P. Colli, M. Fabrizio, G. Gilardi,
Global solution to a singular integrodifferential 
system related to the entropy balance,
{\it Nonlinear Anal.} {\bf 66} (2007) 1949-1979. 

{\pier
\bibitem{BCFG3}
E. Bonetti, P. Colli, M. Fabrizio, G. Gilardi, 
Existence and boundedness of solutions 
for a singular phase field system, 
{\it J. Differential Equations} {\bf 246} (2009) 3260-3295. 
} 

\bibitem{BFR}
E. Bonetti, M. Fr\'emond, E. Rocca,
A new dual approach for a class of phase transitions
with memory: existence and long-time behaviour of solutions, 
{\it J. Math. Pures Appl.} {\bf 88} (2007) 455-481

\bibitem{Brezis}
H. Brezis,
``Op\'erateurs maximaux monotones et semi-groupes de contractions
dans les espaces de Hilbert'',
North-Holland Math. Stud.
{\bf 5},
North-Holland, Amsterdam, 1973.

\bibitem{CC1}
G. Canevari, P. Colli, 
Solvability and asymptotic analysis
of a generalization of the Caginalp phase field system, 
{\it Commun. Pure Appl. Anal.} {\bf 11} (2012) 1959-1982.

\bibitem{CC2}
G. Canevari, P. Colli,  
Convergence properties
for a generalization of the Caginalp phase field system, 
{\it Asymptot. Anal.} {\bf 82} (2013) 139-162.

\bibitem{Fremond}
M. Fr\'emond,
``Non-smooth Thermomechanics'',
Springer-Verlag, Berlin, 2002.

\bibitem{GR}
G. Gilardi, E. Rocca,
Convergence of phase field to phase relaxation governed by the entropy balance with memory,
{\it Math. Methods Appl. Sci.} {\bf 29} (2006) 2149-2179. 

\bibitem{GreenNaghdi}
A.E. Green, P.M. Naghdi, 
A re-examination of the basic postulates of thermo-mechanics, 
{\it Proc. Roy. Soc. Lond. A} {\bf 432} (1991) 171–194.

\bibitem{GriLonSta}
G. Gripenberg, S-O. Londen, O. Staffans,
``Volterra integral and functional equations'',
Encyclopedia Math. Appl.,
{\bf 34},
Cambridge University Press, Cambridge, 1990.

\bibitem{GPmodel}
 M.E. Gurtin, A.C. Pipkin,
{A general theory of heat conduction with finite wave speeds},
{\it Arch. Rational Mech. Anal.}
{\bf 31} (1968) 113-126


\bibitem{Podio1}
P. Podio-Guidugli, 
A virtual power format for thermomechanics, 
{\pier {\it Contin. Mech. Thermodyn.}} {\bf 20} (2009) 479-487.

\bibitem{Simon}
J. Simon,
Compact sets in the space $L^p(0,T; B)$,
{\it Ann. Mat. Pura Appl.~(4)} {\bf 146} (1987) 65-96.

\End{thebibliography}

\End{document}

\bye